\documentclass{article}
\usepackage{amsmath}
\usepackage{amssymb}
\newtheorem{thm}{Theorem}[section]
\newtheorem{lem}[thm]{Lemma}
\newtheorem{prop}[thm]{Proposition}

\title{Associated graded rings of one-dimensional analytically irreducible 
rings II}
\small{\author{ Valentina Barucci\\ Dipartimento di matematica\\
Universit\`a di Roma 1\\ Piazzale Aldo Moro 2\\ 00185 Roma, 
Italy\\
{\em email} barucci@mat.uniroma1.it\\ Ralf Fr\"oberg\\
Matematiska Institutionen\\ Stockholms Universitet\\ 10691 Stockholm, Sweden\\ {\em email} ralff@math.su.se}}

\date{ }

\begin{document}
\maketitle

\begin{abstract}
Lance Bryant noticed in his thesis \cite{lb}, that there was a flaw in our
paper \cite{bf}. It can be fixed by adding a condition, called the BF condition
in \cite{lb}. We discuss some equivalent conditions, and show that they are 
fulfilled for some classes of rings, in particular for our motivating example
of semigroup rings. Furthermore we discuss the connection to a similar
result, stated in more generality, by Cortadella-Zarzuela in \cite{cz}.
Finally we use our result to conclude when a semigroup ring in embedding
dimension at most three has an associated graded which is a complete intersection.

\medskip

\end{abstract}

 \medskip

%%%%%%%%%%%%%%%%%%%%%%%%%%%%%%%%%%%%%%%%%%%%%%%%%%%%%%%%%%%%%%%%%%%%%%%%%%%%
%%%%%%%%%%%%%%%%%%%

%%%%%%%%%%%%%%%%%%%%%%%%%%%%%%%%%%%%%%%%%%%%%%%%%%%%%%%%%%%%%%%%%%%%%%%%%%%%
%%%%%%%%%%%%%%%%%%%%%
\section{The BF condition}
Let $(R,m)$ be an equicharacteristic analytically irreducible and residually 
rational local 1-dimensional domain of embedding dimension $\nu$, multiplicity 
$e$ and residue field $k$. For the problems we study we may, and will, without 
loss of generality suppose that $R$ is complete. So our hypotheses are 
equivalent to supposing $R$ is a subring of $k[[t]]$ with $(R:k[[t]]) \neq 0$. 
Since $k[[t]]$, the integral closure of $R$, is a DVR, every
nonzero element of $R$ has a value, and we let $S=v(R)=\{ v(r); r\in R,r\ne0\}$.
  We denote by $w_0,\ldots,w_{e-1}$  the Apery set of $v(R)$ with respect to 
$e$, 
  i.e., the set of smallest values in $v(R)$ in each congruence class 
$\pmod e$, and we assume $w_j \equiv j \pmod e$.
 
 If $x \in R$ is an element of smallest positive value, i.e. $v(x)=e$, then 
$xR$ is a minimal reduction of the maximal ideal, i.e. $m^{n+1}=xm^n$, for 
$n>>0$. Conversely each minimal reduction of the maximal ideal is a principal 
ideal generated by an element $x$ of value $e$. The smallest integer $n$ such 
that $m^{n+1}=xm^n$ is called the reduction number and we denote it by $r$. 
 
 Observe that, if $v(x)=e$, then Ap$_e(S)=S \setminus (e+S)= v(R) \setminus 
v(xR)$, therefore $w_j \notin v(xR)$, for $j=0, \dots, e-1$.
 
 Consider the $m$-adic filtration $m \supset m^2 \supset m^3 \supset \dots$. 
If $a \in R$, we set ord$(a):= {\rm max}\{i\ |\ a \in m^i\}$. If $s \in S$, we 
consider the semigroup filtration $v(m) \supset v(m^2) \supset \dots$ and set 
vord$(s):= {\rm max}\{i\ |\ s \in v(m^i)\}$. If $a \in m^i$, then $v(a) \in 
v(m^i)$ and so ord$(a) \leq {\rm vord}(v(a))$.
 
According to \cite{lb}, we say that the $m$-adic filtration is {\em 
essentially divisible with respect to the minimal reduction} $xR$   if, 
whenever  $u\in v(xR)$, then
there is an $a\in xR $ with $v(a)=u$ and ord$(a)={\rm vord}(u)$. The $m$-adic 
filtration is {\em essentially divisible} if there
exists a minimal reduction $xR$ such that it is essentially divisible with 
respect to  $xR$.

 We fix for all the paper the following notation. Set, for $j=0, \dots, e-1$,
$b_j=\max\{ i| w_j\in v(m^i)\}$, and let $c_j=\max\{ i|w_j\in v(m^i+xR)\}$. 
Note that the numbers $b_j$'s  do not depend on the minimal reduction $xR$, 
on the contrary the $c_j$'s depend on $xR$.

%\begin{lem} Suppose that a set $V=\{v_1, \dots,v_n\} $ is the disjoint union of  subsets $A\cup B$ and $C \cup D$, with $|A|=|D|=h$ and $|B|=|C|=n-h$. Then $|A \cap C|=|B \cap D|$.\end{lem}\noindent{\bf Proof.} Suppose that $A \cap C=\{v_1, \dots,  v_s\}$, then $v_1, \dots,  v_s \notin D$ and $v_{s+1}, \dots,  v_h \in D$. So $D$ has $h-(h-s)=s$ elements from $B$. It follows that $|A \cap C|\leq|B \cap D|$. By simmetry we get the equality. $\Box$

\begin{lem} \label{stupid}
If $I$ and $J$ are ideals of $R$, then $v(I+J)=v(I)\cup v(J)$
is equivalent to $v(I\cap J)=v(I)\cap v(J)$.
\end{lem}
\noindent {\bf Proof.} Let $V=v(I+J)\setminus v(I\cap J)$. Then 
$$V=(v(I)\setminus v(I\cap J))\cup(v(I+J)\setminus v(I))=
(v(J)\setminus v(I\cap J))\cup(v(I+J)\setminus v(J))$$
 and both unions are
disjoint. Since $(I+J)/J\simeq I/I\cap J$, 
we get that
$|v(I+J)\setminus v(J)|=|v(I)\setminus v(I\cap J)|$ and similarly that $|v(I+J)\setminus v(I)|=|v(J)\setminus v(I\cap J)|$. Suppose that $v(I\cap J)\subsetneq v(I)\cap v(J)$, i.e. that there is a value $v_0 \in (v(I)\setminus v(I\cap J)) \cap (v(J)\setminus v(I\cap J))$. Thus $v_0 \notin (v(I+J)\setminus v(J))$ and by cardinality reasons also $(v(I+J)\setminus v(I)) \cap (v(I+J)\setminus v(J)) \neq \emptyset$, i.e. $v(I+J) \supsetneq v(I)\cup v(J)$. The other implication is symmetric and we get the claim. $\Box$

\begin{prop}\label{essdiv} Let $xR$ be a minimal reduction of $m$.
Then the  following conditions are equivalent:

{\rm (1)} The $m$-adic filtration is essentially divisible with respect to $xR$.

{\rm (2)} $v(m^i\cap xR)= v(m^i)\cap v(xR)$, for all $i\geq 0$.

{\rm (3)} $v(m^i+xR)=v(m^i)\cup v(xR)$ for all $i\geq 0$.

{\rm (4)}  $b_j=c_j$ for $j=0,\ldots,e-1$.

\end{prop}
\noindent {\bf Proof.} (1)$\Rightarrow$(2): Let $i \geq 0$ and $u \in v(m^i) 
\cap v(xR)$. Then $u \in v(xR)$ and vord$(u) \geq i$. By (1) there exists 
$a\in xR $ with $v(a)=u$ and ord$(a)={\rm vord}(u)$. Thus $a \in m^i \cap xR$ 
and so  $v(m^i\cap xR) \supseteq v(m^i)\cap v(xR)$. Since the other inclusion 
is trivial, we get an equality.

(2)$\Rightarrow$(1): If $u \in v(xR)$ and vord$(u)= i$, then $u \in  v(m^i)
\cap v(xR)$, and by (2), $u \in v(m^i\cap xR)$. So there is $a \in m^i \cap 
xR$ with $v(a)=u$. For such $a$, $i \leq {\rm ord}(a) \leq {\rm vord }(u)=i$, 
and so  ord$(a)=i$.

That (2) and (3) are equivalent follows from Lemma \ref{stupid} with
$I=m^i$ and $J=xR$.

(3)$\Rightarrow$(4): Since $m^i\subseteq m^i+xR$, we have 
$v(m^i)\subseteq v(m^i+xR)$, so $b_j\le c_j$. Suppose that $b_j<c_j$ for some
$j$. Then $w_j\in v(m^{c_j}+xR)\setminus v(m^{c_j})$. Since $w_j\notin v(xR)$,
we get that $v(m^{c_j})\cup v(xR)$ is strictly included in $v(m^{c_j}+xR)$.\\
(4)$\Rightarrow$(3): If $u\in v(m^i+xR)\setminus v(xR)$, then 
$u\in v(R)\setminus v(xR)={\rm Ap}_ev(R)$, so $u=w_j$ for some $j$. Then
$w_j\in v(m^i+xR)\setminus v(m^i)$, so $b_j<c_j$. $\Box$
 
 \bigskip
 
 Observe that if $R=k[[t^{n_1}, \dots, t^{n_{\nu}}]]$ is a semigroup $k$-algebra 
and $I$, $J$ are ideals generated by monomials, then    $v(I\cap J)=v(I)\cap 
v(J)$ (and $v(I+J)=v(I)\cup v(J)$).   This follows from the fact that if
$I=(t^{i_1},\ldots,t^{i_k})$ is generated by monomials, then $v(I)=
\langle i_1,\ldots,i_k\rangle$. So, if we choose for the maximal ideal of 
$R$ a monomial minimal reduction, by Proposition \ref{essdiv} we have that the 
$m$-adic filtration  is essentially divisible with respect to such a 
reduction. If we choose a different minimal reduction this is not always the 
case, as the following example shows.
 
 \bigskip
 \noindent{\bf Example} Let $R=k[[t^6,t^7,t^{15}]]$. By what we observed above, 
the $m$-adic filtration is essentially divisible with respect to the minimal 
reduction $t^6R$. On the contrary, it is not essentially divisible with 
respect to the minimal reduction $(t^6+t^7)R$, because $v(m^3+ (t^6+t^7)R) 
\nsubseteq  v(m^3) \cup v( (t^6+t^7)R)$ and we can apply Proposition 
\ref{essdiv} (3). As a matter of fact, $t^{21}-(t^6+t^7)t^{15} \in m^3+ 
(t^6+t^7)R$, thus $22 \in v(m^3+ (t^6+t^7)R)$, but $22 \notin v(m^3) \cup v( 
(t^6+t^7)R)$.
 
 This example shows also that the numbers $c_j$'s depend on the minimal 
reduction. Considering $w_4=22$, with respect to the minimal reduction $t^6R$, 
we get $b_4=c_4=2$, but with respect to $(t^6+t^7)R$, we get $2=b_4<c_4=3$.
 
 \bigskip
 
 In \cite {bf}, we called a set 
$f_0,\ldots,f_{e-1}$ of elements of $R$ an {\em Apery basis} if  $v(f_j)
\equiv j\pmod e$ and ord$(f_j)=b_j$, for all $j$, $j=0, \dots, e-1$ and 
claimed that for all $i\ge0$, $m^i$ is a free $W$-module generated by elements 
of the 
form $x^{h_j}f_j$, where $xR$ is a minimal reduction of $m$ and $W=k[[x]]$. In 
\cite{lb} Lance Bryant showed that this is not always true, considering the 
example  
$R=k[[t^6,t^8+t^9,t^{19}]]$ with char$(k)=0$. 
Here  $e=6$  and  $v(R)$ has Apery set $0,8,16,19,27,29$.    Setting: 
$ x=t^6, W=k[[t^6]]$ and 
$f_0=1,  f_1= t^8+t^9, f_2= t^{16}+2t^{17}+t^{18}, f_3=t^{19},f_4=t^{27}+t^{28}, 
f_5=t^{29}$
 he gets $m^3= x^3f_0W+x^2f_1W+xf_2W+ gW+xf_4W+xf_5W$
where $g= (t^8+t^9)^3-(t^6)^4 =3t^{25}+3t^{26}+t^{27}\in m^3$. On the other hand 
$x^hf_3=t^6t^{19}=t^{25} \in m^2 \setminus m^3$.

According to \cite{lb}, we say that the $m$-adic filtration satisfies the   BF 
{\em condition} if there exists a minimal reduction $xR$ of $m$ and a set of 
elements $\{f_0, \dots, f_{e-1}\}$ of $R$ with $v(f_j)=w_j$  such that  each 
power of  $m$ is a free $k[[x]]$-module generated by elements of the form 
$x^{h_j}f_j$. 
 
 The BF condition depends on the choice of the  elements $\{f_0, \dots, 
f_{e-1}\}$ and on the reduction. In \cite{bf}  we noted that, if $R=k[[t^4,
t^6+t^7,t^{13}]] $, with char$(k) \neq2$,  then 
 ${\rm Ap}_4(v(R))=\{0,6,13,15\}$ and 
setting $f_0=1, \  f_1= t^6+t^7, \ f_2= 2t^{13}+t^{14},\  f_3=t^{15}$,
$ x=t^4$, $W=k[[t^4]]$,
  we get that each power of the maximal ideal is a free $W$-module generated 
by elements of the form 
$x^{h_j}f_j$. For example:  
  $$m= xf_0W+f_1W+f_2W+f_3W$$
  $$m^2= x^2f_0W+xf_1W+f_2W+xf_3W$$
$$m^3=xm^2= x^3f_0W+x^2f_1W+xf_2W+xf_3W$$
 If we replace $f_2$ with $t^{13}$, since $t^{13} \in m \setminus m^2$, we don't 
have the free basis of the requested form for $m^2$. Thus this example shows 
that the BF condition depends on the choice of the  elements $\{f_0, \dots, 
f_{e-1}\}$.
To show that the BF condition  
 depends on the reduction, we can consider the example above, $R=k[[t^6,t^7,
t^{15}]]$.  We get that $f_0=0,f_1=t^7,f_2=t^{14},f_3=t^{15},f_4=t^{22},
f_5=t^{29}$
is an Apery basis but, choosing the minimal reduction $xR=(t^6+t^7)R$, $m^4$ 
is not a free $k[[x]]$-module generated by elements of the form $x^{h_j}f_j$, 
because Ap$_6(v(m^4))= \{24,25,26,27,28,35\}$ and an element of the form 
$x^{h_j}f_j$  of value $28$ is $(t^6+t^7)t^{22}$, which   is not in $m^4$.

\begin{prop} \label{bf} Let $W=k[[x]]$, where $xR$ is a minimal reduction of 
$m$ and let 
$f_0,\ldots,f_{e-1}$ be elements of $R$ with $v(f_j)\equiv j\pmod e$.
Then the following conditions are equivalent:

\noindent {\rm (1)} For all $i\ge0$, $m^i$ is a free $W$-module generated by 
elements of the 
form $x^{h_j}f_j$.\\
{\rm (2)} For all $i\ge0$, {\rm Ap}$_e(v(m^i))=\{v(x^{h_j}f_j)\}$ for some 
$x^{h_j}f_j\in m^i$, $j=0,\ldots,e-1$.\\
{\rm (3)} If $\sum_{j=0}^{e-1}d_j(x)f_j\in m^i$ with $d_j(x)\in W$ for all $j$, 
then
$d_j(x)f_j\in m^i$ for each $j$.
\end{prop}

\noindent {\bf Proof.} (1)$\Rightarrow$(3): Let $
a=\sum_{j=0}^{e-1}d_j(x)f_j\in m^i$.
Since $\{ x^{h_j}f_j\}$ is a free basis for $m^i$, we also have 
$a=\sum_{j=0}^{e-1} d'_j(x)x^{h_j}f_j$ for some $d'_j(x)$, and 
$d_j(x)=d'_j(x)x^{h_j}$. Now $x^{h_j}f_j\in m^i$, so $d_j(x)f_j\in m^i$.\\
(3)$\Rightarrow$(2): Let $u\in {\rm Ap}_e(v(m^i))$, so $u=v(a)$ for some
$a\in m^i$. We have $a=\sum_{j=0}^{e-1}d_j(x)f_j$, with $d_j(x)f_j\in m^i$
for all $j$. Let $v(a)\equiv v(f_j)\pmod e$. Then $v(a)=v(d_j(x)f_j)$. Let
$d_j(x)=\sum_{i\ge l} k_ix^i$, with $k_i\in k, k_l\ne0$. Then we claim that
ord$(d_j(x)f_j)={\rm ord}(x^lf_j)$. Suppose that 
$x^lf_j\in m^h\setminus m^{h+1}$. Then $d_j(x)f_j\in m^h$ since all summands do.
If $d_j(x)f_j\in m^{h+1}$, then 
$k_lx^lf_j=d_j(x)f_j-\sum_{i\ge l+1}k_ix^if_j\in m^{h+1}$, a contradiction.
Thus $v(a)=v(x^lf_j)$, $x^lf_j\in m^i$.\\
(2)$\Rightarrow$(1):  By Lemma 2.1 (1) of \cite{bf}. $\Box$

\begin{prop} \label{bfstronger} If the $m$-adic filtration satisfies the BF 
condition, it is essentially divisible.
\end{prop}

\noindent {\bf Proof.} Let $xR$ be a minimal reduction of $m$ and let
$f_0,\ldots,f_{e-1}$ be elements in $R$ satisfying the BF condition, i.e.  
condition (2)
in Proposition \ref{bf}.  We claim that condition (2) in Proposition \ref 
{essdiv} is satisfied.
 Let $v\in v(m^i)\cap v(xR)$, $v=v_j+le$, with 
$v_j\in{\rm Ap}_e(v(m^i))$, for some $l \geq 0$.   We have
$v_j=v(x^{h_j}f_j)$, for some $j$. Thus $x^{h_j+l}f_j\in m^i\cap xR$ and
$v(x^{h_j+l}f_j)=v$.   Note that $h_j+l >0$.  $\Box$

\bigskip

There are several cases in which the BF condition holds.

\begin{prop}\label{cases}
The BF-condition holds for the $m$-adic filtration in each of the following 
cases:

{\rm (1)} $R$ is a semigroup $k$-algebra.

{\rm (2)} The reduction number $r$   is at
most 2.

{\rm (3)}  The embedding dimension $\nu$ is at most 2.
\end{prop}
\noindent {\bf Proof.} (1): Let $R=k[[t^{n_1}, \dots, t^{n_ \nu}]]$ and 
Ap$(v(R))=\{w_0, \dots, w_{e-1}\}$. Choosing the monomial Apery basis 
$f_j=t^{w_j}$, for $j=0, \dots, e-1$ and the monomial minimal reduction 
$xR=t^{n_1}R=t^eR$, if Ap$(v(m^i))=\{w_0+h_0e, \dots, w_{e-1}+h_{e-1}e\}$, then 
$m^i$ is a free $k[[t^e]]$-module generated by $t^{eh_j}f_j= t^{h_je+w_j}$.
 
\noindent
(2): Let  $xR$ is a minimal reduction of 
$m$ and let $f_0,\ldots,f_{e-1}$ be an Apery basis of $R$. Then the Apery sets 
of $v(m^i))$, with $i \leq 2$ can always be realized as in Proposition 
\ref{bf} (2).  In fact, for $v(m^2)$, note that $v(x^2f_0)=2e \in 
{\rm Ap}(v(m^2))$. Moreover, if $f_j \in m \setminus m^2$, then $v(xf_j) 
\in {\rm Ap}(v(m^2))$ and if $f_j \in m^2$, then $v(f_j) \in {\rm Ap}(v(m^2))$. 
If $i\ge2$, then $m^{i+1}=xm^i$, which gives the claim. 

\noindent
(3)
In the plane case,
setting $m=\langle x,y \rangle$, using the Weierstrass Preparation
Theorem, we noted in \cite[Section 2]{plane} that $R$  is a $W$-module generated  by
$1,y,y^2,...,\\ y^{e-1}$ and replacing each $y^j$ with a suitable $y_j=y^j+ \phi(x,y)$ ($\phi(x,y)\in m^j$), we get an Apery basis for $R$. Consider a power $m^i$ of the maximal ideal.  Using the above observation,   $m^i$ is generated as $W$-module by $ x^i,x^{i-1}y,x^{i-2}y^2,\dots, y^i, y^{i+1},\\
\dots,y^{i(e-1)}$.
Now working on the powers $y^j$ as we do in \cite{plane}, we can modify the generators, getting the $e$ elements  $x^i,x^{i-1}y,x^{i-2}y_2,\dots, y_{e-1}$, which are still in $m^i$, are of the requested form and such that their values form an Apery set for $v(m^i)$.
$\Box$ 

\bigskip
\noindent {\bf Example} Consider $R=\mathbb C[[t^6,t^8+t^9]]$. Setting $x=t^6$, $y=t^8+t^9$, as in \cite{plane}, we can see that an Apery basis for $R$ is $1,y,y_2=y^2,y_3=y^3-x^4= 3t^{25}+..., y_4=y^4-x^4y=5t^{33}+...,y_5=y^5-x^4y^2=
5t^{41}+...$.
Considering for example $m^3$, we see it is a free $W$-module generated by
$x^3, x^2y, xy_2, y_3, y_4,y_5$.

 %%%%%%%%%%%%%%%%%%%%%%%%%%%%%%%%%%%%%%%%%%%%%%%%%%
\section{The associated graded ring}
Let gr$(R)$ be  the associated graded ring with respect to the $m$-adic 
filtration, gr$(R)= \bigoplus_{i \geq 0}m^i/m^{i+1}$.  The CM-ness of gr$(R)$ is 
equivalent to the existence of a nonzerodivisor in the homogeneous maximal 
ideal. If such a nonzerodivisor exists, then $x^*$, the image of $x$ in 
gr$(R)$ (where $x$ is any element of value $e$) is a nonzerodivisor. We fix this notation and denote by Hilb$_R(z)= \sum_{i \geq 0}l_R(m^i/m^{i+1})z^i$ the Hilbert series of $R$ and by Hilb$_{R/xR}(z)= \sum_{i \geq 0}l_R(m^i+xR/m^{i+1}+xR)z^i$ the Hilbert series of $R/xR$. Recall that $$(1-z){\rm Hilb}_R(z) \leq {\rm Hilb}_{R/xR}(z)$$and the equality holds if and only if gr$(R)$ is CM (cf. e.g. 
[3] or [4]).

We start noting that, if gr$(R)$ is CM, then the conditions analyzed in the 
previous section are equivalent. 

\begin{prop} \label{equivunderCM} If {\rm gr}$(R)$ is CM, then the $m$-adic 
filtration is essentially divisible
 if and only if it satisfies the BF condition.
\end{prop}

\noindent{\bf Proof.} 
Suppose that the $m$-adic filtraion is essentially divisible
with respect to $xR$. We claim that  there exist $f_0,\ldots,f_{e-1}$ in $R$
satifying  condition (2)  of Proposition \ref{bf}.
 If $n\ge r$, where $r$ is the reduction number, then 
$m^n\subseteq xR$. Thus, if $u\in{\rm Ap}_e(v(m^n))$, $u\equiv j\pmod{e}$,
then there exist $a\in R$, $a=xa'$, with $v(a)=u$ and ord$(a)=n$. We have
$v(a')=u-e$ and ord$(a')={\rm ord}(a)-1$, because gr$(R)$ is CM. Now there 
are two possibilities. If $v(a')\notin v(xR)$, i.e. $v(a')=w_j$, we
choose $f_j=a'$. If $v(a')\in v(xR)$, then, since $R$ is essentially
divisible, there exist $b\in xR$, $b=xb'$, with $v(b)=v(a')$ and ord$(b)=
{\rm ord}(a')$. Moreover $b\in{\rm Ap}(v(m^{n-1}))$, because otherwise
$u-2e\in v(m^{n-1})$ and $u-e\in v(m^n)$, a contradiction. Continuing in
this way we arrive to get the element $f_j$ requested.

\bigskip

We denote by $R'$ the first neighborhood ring or the blowup of $R$, i.e. the 
overring $\bigcup_{n \geq0}(m^n:m^n)$. It is well known that, if $v(x)=e$, 
$R'=R[x^{-1}m]= \bigcup_{i \geq 0}\{yx^{-i}; y \in m^i\}$, cf. \cite{l}.
Let $w'_0, \dots, w'_{e-1}$ be the Apery set of $v(R')$ with respect to $e$, 
with $w'_j \equiv j \pmod e$. For each $j$, $j=0, \dots, e-1$, define as in 
\cite {bf} $a_j$ by $ w'_j=w_j-a_je$.

If $f_j \in m^i$, then $f_jx^{-i} \in R'$, so $v(f_jx^{-i})=w_j-ie \in v(R')$. 
It follows that $w_j-b_je \in v(R')$. Since $ w'_j=w_j-a_je$ is the smallest in 
$v(R')$, in its congruence class $\pmod e$, we have that $a_j \geq b_j$, for   
$j=0, \dots, e-1$.

In \cite [Theorem 2.6]{bf} we stated the following:
 The ring {\rm gr}$(R)$ is CM if and only if $a_j =b_j$, for   
$j=0, \dots, e-1$.

 As Lance Bryant pointed out, the proof of that theorem given in \cite{bf} 
works under the assumption that the $m$-adic filtration satisfies the BF 
condition. 

 \begin{thm}\label{thm}
 If $R$ satifies the BF condition then {\rm gr}$(R)$ is CM if and only if $a_j =b_j$, for  $j=0, \dots, e-1$.

\end{thm}

\noindent {\bf Proof.} If the BF condition is satisfied, the proof given in 
\cite{bf} holds.

\bigskip

In \cite {cz} T.~Cortadellas  and  S.~Zarzuela proved, in more general 
hypotheses for $R$, a criterion for the CM-ness of gr$(R)$. They consider the 
microinvariants of J.~Elias, i.e. the numbers $\epsilon_j$ which appear in the 
decomposition of the torsion module 
$$R'/R= \bigoplus_{j=0}^{e-1} W/x^{\epsilon_j}W$$ 
where $R'$ is the blowup, $xR$ a minimal reduction of $m$ and $W=k[[x]]$. With 
our hypotheses and notation, they show in particular that gr$(R)$ is CM if and 
only if $c_j=\epsilon_j$, for $j=0, \dots, e-1$,  \cite[Theorem 4.2]{cz}.
Comparing their result with ours, we see that they are coherent but different. 
In fact, if the $m$-adic filtration satisfies the BF condition, then,  for 
$j=0, \dots, e-1$, $\epsilon_j=a_j$ by \cite[Proposition 2.5]{bf} and $b_j=c_j$ 
by Propositions \ref{essdiv} and \ref{bfstronger}, so their result coincide 
with ours.  The hypotheses on the ring in their result are more general, but 
the numbers $c_j$'s depend on the minimal reduction. On the other hand, the 
numbers $a_j$'s and $b_j$'s which we consider do not depend on the minimal 
reduction and in our criterion the CM-ness of gr$(R)$ can be read off just 
looking at the semigroup filtration $v(m^0) \supset v(m) \supset v(m^2) \supset 
\dots$. As a matter of fact, since $R'= x^{-n}m^n$, for $n>>0$, $v(R') =
v(m^n)-ne$, for $n>>0$, so the $a_j$'s which relate the Apery sets of $v(R)$ 
and $v(R')$, can be read in the semigroup filtration $\{v(m^i)\}_{i \geq 0}$. 

\bigskip

We give now some applications. Given an analytically irreducible ring 
satisfying our hypotheses, we denote by $a_j(R)$ and $b_j(R)$ the numbers 
defined above. 
\begin{prop} Let $R$ and $T$ be rings satifying the BF condition, with the 
same multiplicity $e$ and  with
$a_j(R)=a_j(T)$, $b_j(R)=b_j(T)$, for $j=0, \dots, e-1$. If {\rm gr}$(R)$ is CM,
then also gr$(T)$ is CM and $R$ and $T$ have the same Hilbert series.

\end{prop}

\noindent {\bf Proof.} Since gr$(R)$ is CM, by Theorem \ref{thm},  
$a_j(R)=b_j(R)$, for $j=0, \dots, e-1$. So also $a_j(T)=b_j(T)$, for $j=0, 
\dots, e-1$ and gr$(T)$ is CM.  If $xR$ 
(respectively $yT$) is a minimal reduction of the maximal ideal of $R$ 
(respectively of $T$), then, since $b_j(R)=c_j(R)$ and $b_j(T)=c_j(T)$ (cf. 
Proposition \ref{essdiv}), the Hilbert series of $R/xR$ and $T/yT$ are the 
same.  Since 
Hilb$_{R/xR}(z)=(1-z){\rm Hilb}_R(z)$ and Hilb$_{T/yT}(z)=(1-z){\rm Hilb}_R(z)$,
also the Hilbert series of $R$ and $T$ are the same. $\Box$

\bigskip

Sometimes we can use the BF condition to draw conclusions about when
gr$(R)$ is a complete intersection (CI). We will use that if $x\in R $ is
a nonzerodivisor in $R$ such that $x^*$ is a nonzerodivisor in gr$(R)$, then
gr$(R/xR)={\rm gr} (R)/(x^*)$, \cite[Lemma(b)]{he2}.

\bigskip

\noindent
{\bf Example}
 If $R=k[[X,Y]]/(f)$ is a plane branch, then gr$(R)=k[X,Y]/(f^*)$, where
$f^*$ is the image of $f$ in gr$(R)$, so gr$(R)$ is a complete intersection.
The semigroups $S$ for which $k[[S]]$ is a CI were determined in \cite{de}.
If gr$(k[[S]])$ is a CI, then necessarily $k[[S]]$ 
is a CI \cite[Corollary 2.4]{va}.
If $S$ is generated by three elements and is a CI, the generators are of 
the form
$na,nb,n_1a+n_2b$, $a<b$, \cite{he1} or (with an easier proof) 
\cite[Lemma 1]{wa}. Then 
$$k[[S]]=k[[X,Y,Z]]/(X^b-Y^a,Z^n-X^{n_1}Y^{n_2})$$ It is determined in \cite{he2}
when gr$_m(k[[S]])$ is a CI when $S$ is 3-generated. The result is

\noindent
a) $S=\langle na,nb,n_1a\rangle$. 

\noindent
b) $S=\langle na,nb,n_1a+n_2b\rangle$, $na<n_1a+n_2b<nb$, $n\le n_1+n_2$.

\noindent
c) $S=\langle na,nb,n_1a+n_2b\rangle$, $na<nb<n_1a+n_2b$, $n\le n_1+n_2$.

\smallskip\noindent
Let $x=t^{na},y=t^{nb},z=t^{n_1a+n_2b}$.

\smallskip\noindent
In case a), if $n<n_1$, 
gr$(k[[S]]/(x))\cong k[Y,Z]/(Y^a,Z^n)$.
An Apery basis for $k[[S]]$ is $\{ y^iz^j;0\le i<a,0\le j<n\}$. Suppose
  $R=k[[t^{na},g_2,g_3]]$ with $v(g_2)= nb,v(g_3)=n_1a$, and that
$\{ g_2^ig_3^j;0\le i<a,0\le j<n\}$ is an Apery basis for $R$, and that $R$
satisfies the BF condition. Then $x=t^{na}$ is a minimal reduction also of
the maximal ideal of $R$, and the $a_j$'s and $b_j$'s are the same for
$k[[S]]$ and $R$, so gr$(R)$ is CM, and in particular $x^*$ is a
nonzerodivisor in  gr$(R)$. We have that gr$(R)$ is a CI if and only if
gr$(R/xR)={\rm gr}(R)/(x^*)$ is a CI.
Since $v(g_2^ig_3^j)\notin v(xR)$ if $0\le i<a,0\le j<n$, and they all have
values in different congruence classes $\pmod{v(x)}$, we get that 
 ${\rm gr}(R)/(x^*) \cong  {\rm gr}(k[[S]])/(x^*) \cong  k[Y,Z]/(Y^a,Z^n)$.
  Thus gr$(R)$ is a CI.   A concrete
example is   $R=k[[t^6,t^8+ct^{13}+dt^{19},t^9]]$, $c,d \in k$.

If $n_1<n$, then gr$(k[[S]]/(z))=k[X,Y]/(Y^a,X^{n_1})$, and $\{ y^ix^j;0\le i<a,0\le
j<n_1\}$ is an Apery basis for $k[[S]]$. Suppose
$R=k[[t^{n_1a},g_2,g_3]]$ with $v(g_2)=na,v(g_3)=nb$, and that
$\{ g_3^ig_2^j;0\le i<a,0\le j<n_1\}$ is an Apery basis for $R$,  and  that $R$
satisfies the BF condition. As above we get that gr$(R)$ is a CI. A concrete 
example is $k[[t^6,t^9+ct^{11},t^4]]$, $c \in k$.

\smallskip\noindent
In case b) an Apery set is $\{ y^iz^j;0\le i<a,0\le j<n\}$. Suppose 
$R=k[[t^{na},g_2,g_3]]$, $v(g_2)=n_1a+n_2b,v(g_3)=nb$, and that 
$\{ g_3^ig_2^j;0\le i<a,0\le j<n\}$ is an
Apery set for $R$,  and that   $R$ satifies the BF condition. Reasoning as above, we
get that gr$(R)$ is a CI. A concrete example is 
 $k[[t^6,t^7+ct^{11},t^9]]$, $c \in k$.

\smallskip\noindent
In case c) an Apery set is $\{ y^iz^j;0\le i<a,0\le j<n\}$. Suppose 
$R=k[[t^{na},g_2,g_3]]$, $v(g_2)=nb,v(g_3)=n_1a+n_2b$, and that 
$\{ g_2^ig_3^j;0\le i<a,0\le j<n\}$ is an
Apery set for $R$,  and that $R$ satifies the BF condition. Reasoning as above, we
get that gr$(R)$ is a CI. A concrete example is  
$k[[t^4,t^6,t^7+ct^9]]$, $c \in k$.

\medskip

We end with some questions:

\begin {enumerate}
\item Does the converse of Proposition \ref{bfstronger} hold?
\item Is Theorem \ref{thm} true, without assuming the BF-condition?
\item Is always $\epsilon_j=a_j$, for $j=0,\dots,e-1$ without assuming the BF-condition?

\end{enumerate}

\end{document}